\documentclass[a4paper,11pt]{article}
\usepackage{latexsym}
\newtheorem{thm}{Theorem}[section]
\newtheorem{lem}[thm]{Lemma}

\newtheorem{cor}[thm]{Corollary}
\newtheorem{prop}[thm]{Proposition}

\newtheorem{rem}[thm]{Remark}

\input epsf

\title{On Minkowski's bound for lattice-packings}
\author{Roland Bacher\footnote{Support 
from the Swiss National Science Foundation is
gratefully acknowledged.}}
\date{}
\begin{document}
\maketitle

%\par fichier dense3.tex dans reseaux, derni\`ere modification: aout 04 \bigskip
\par {\it Abstract: We give a new proof of the Minkowski-Hlawka bound
on the existence of dense lattices.
The proof is based on an elementary method for constructing
dense lattices which is almost effective.}
\footnote{Math. class.: 10E05, 10E20. Keywords: Lattice packing}
 
\section{Introduction}
                          
Let $\mu\geq 2$ be a strictly positive integer. A {\it $\mu-$sequence} is a
sequence $s_0=1,s_1,s_2,\dots$ of strictly positive integers such that 
the $n-$dimensional lattice
$$\Lambda_n=\{(z_0,z_1,\dots,z_n)\in{\mathbf Z}^{n+1}\ \vert \ 
\sum_{k=0}^n s_kz_k=0\}=(s_0,\dots,s_n)\cap{\mathbf Z}^{n+1}$$
has minimum $\geq \mu$ for all $n\geq 1$.
Since $\hbox{det}(\Lambda_n)=\sum_{k=0}^n s_k^2$ we get a lower bound for 
the center-density 
$$\delta(\Lambda_n)=\sqrt{\frac{(\hbox{min}\ \Lambda_n)^n}{4^n\ \hbox{det}
\ \Lambda_n}}\geq \sqrt{\frac{\mu^n}{4^n\ \sum_{k=0}^n s_k^2}}$$
(or for the density $\Delta(\Lambda_n)=\delta(\Lambda_n)\pi^{n/2}/(n/2)!$)
of the $n-$dimensional lattice $\Lambda_n$ associated to a $\mu-$sequence.

\begin{thm} \label{mainA} Given an integer $\mu\geq 2$ as above there exists a
$\mu-$sequence $s_0=1,s_1,\dots$ satisfying for all $n\geq 1$
$$s_n\leq 1+\sqrt{\mu-2}\sqrt{\mu-1+n/4}^n\frac{\sqrt{\pi}^n}{(n/2)!}
\leq \sqrt{\mu}\sqrt{\mu+n/4}^n\frac{\sqrt{\pi}^n}{(n/2)!}\ .$$
\end{thm}

\begin{rem} \label{rem1} (i) The condition $s_0=1$ is of no real importance
and can be omitted after minor modifications. It is of course 
also possible (but not very useful) to consider sequences
with coefficients in ${\mathbf Z}$.

\ \ (ii) Any subsequence $s_{i_0}=s_0,s_{i_1},s_{i_2},\dots$
of a $\mu-$sequence is again a $\mu-$sequence. Reordering the
terms of a $\mu-$sequence (in increasing order) yields of course
again a $\mu-$sequence. 

\ \ (iii) The lattices associated to a $\mu-$sequence are generally
neither perfect nor eutactic (cf. \cite{M} for definitions) 
and one can thus generally improve their densities by suitable deformations.
\end{rem}

The proof of Theorem \ref{mainA} is very elementary and
consists of an analysis of the  \lq\lq greedy algorithm''
which constructs the first $\mu-$sequence with respect to
the lexicographic order on sequences. An easy analysis shows that the
lexicographically first sequence satisfies the first inequalities
of Theorem \ref{mainA}. The greedy algorithm, although very simple, 
is however quite useless for applications because of astronomical memory
requirements (which can be lowered at the price of an astronomical amount of
computations).

$\mu-$sequences satisfying the inequalities of Theorem \ref{mainA}
yield rather dense lattices as shown by the next result.

\begin{cor} \label{corA} For any $\mu\geq 2$, there exists a $\mu-$sequence
$(s_0,s_1,\dots,s_n)\in{\mathbf Z}^{n+1}$ such that the density of
the associated lattice $\Lambda_n=(s_0,\dots,s_n)^\perp\cap{\mathbf Z}^{n+1}$
satisfies
$$\Delta(\Lambda_n)
\geq \frac{(1+n/(4\mu))^{-n/2}}{2^{n}\sqrt{(n+1)\mu}}\ .$$
\end{cor}

\begin{rem} \label{rem2} Taking $\mu\sim n^2/4$ 
we get the existence of lattices
in dimension $n$ (for large $n$) with density $\Delta$ roughly at
least equal to
$$\frac{1}{2^{n-1}\ n\ \sqrt{(n+1)\ e}}\ .$$

This is already close to the Minkowski-Hlawka bound (which shows the
existence of lattices
with density at least $\zeta(n)\ 2^{1-n}$, cf. formula (14) in  \cite{CS},
Chapter 1. The best known lower bound concerning densities of
lattice packings
(together with a very nice proof) seems to be due to Ball
and asserts the existence of $n-$dimensional lattices with density 
at least $2(n-1)2^{-n}\zeta(n)$, see \cite{Ball}.
\end{rem}

A more careful analysis of $\mu-$sequences yields
the following result.

\begin{thm} \label{mainB} For every $\epsilon>0$, there 
exist $n-$dimensional lattices with density
$$\Delta\geq\frac{1-\epsilon}{2^n\ \sum_{k=1}^\infty e^{-k^2\pi}}\sim
(1-\epsilon)\ 23.1388\ 2^{-n}$$
for all $n$ large enough.
\end{thm}

\section{Definitions}

For the convenience of the reader this section contains all needed facts 
concerning lattices. Reference for lattices and
lattice-packings are \cite{CS} and \cite{M}. 

An {\it $n-$dimensional lattice} is a discret-cocompact subgroup $\Lambda$ 
of the $n-$dimensional Euclidean vector space ${\mathbf E}^n$ (with
scalar product denoted by $\langle\ ,\ \rangle$). The {\it determinant}
of a lattice $\Lambda$ is the square of the volume of a fundamental
domain ${\mathbf E}^n/\Lambda$ and equals
$\det(\langle b_i,b_j\rangle)$ where $b_1,\dots,b_n$ denotes a 
${\mathbf Z}-$base of $\Lambda$. The {\it norm} of a lattice element
$\lambda\in\Lambda$ is defined as $\langle \lambda,\lambda\rangle$
(and is thus the squared Euclidean norm of $\lambda$). A lattice
$\Lambda$ is {\it integral} if the scalar product takes only integral values on 
$\Lambda\times \Lambda$.
The {\it minimum} 
$$\hbox{min}\ \Lambda\ =\hbox{min}_{\lambda\in\Lambda\setminus\{0\}}
\langle \lambda,\lambda\rangle$$
of a lattice $\Lambda$ is the norm
of a shortest non-zero vector in $\Lambda$. The {\it density}
$\Delta(\Lambda)$ and the {\it center-density}
$\delta(\Lambda)$ of an $n-$dimensional lattice $\Lambda$ are defined as
$$\Delta(\Lambda)=\sqrt{\frac{(\hbox{min}\ \Lambda)^n}{4^n\ \hbox{det}
\ \Lambda}}\ V_n\qquad \hbox{and}\qquad 
\delta(\Lambda)=\sqrt{\frac{(\hbox{min}\ \Lambda)^n}{4^n\ \hbox{det}
\ \Lambda}}$$
where $V_n=\pi^{n/2}/(n/2)!$ denotes the volume of the $n-$dimensional
unit-ball in ${\mathbf E}^n$.
These two densities are proportional (for a given dimension $n$) and
measure the quality of the sphere-packing associated to the
lattice $\Lambda$ obtained by packing $n-$dimensional Euclidean balls
of radius $\sqrt{\hbox{min}\ \Lambda/4}$ centered at all points of $\Lambda$.

Given an $n-$dimensional lattice $\Lambda\subset {\mathbf E}^n$ the subset
$$\Lambda^\sharp=\{x\in{\mathbf E}^n\ \vert\ \langle x,\lambda\rangle\in{\mathbf Z}
\quad \forall \lambda\in\Lambda\}$$
is also a lattice called the {\it dual lattice} of $\Lambda$. The lattice
$\Lambda$ is integral if and only if $\Lambda\subset \Lambda^\sharp$. For an
integral lattice the {\it determinant group} $\Lambda^\sharp/\Lambda$ is
a finite abelian group consisting of $(\hbox{det}\ \Lambda)$ elements.

A sublattice $M\subset \Lambda$ is {\it saturated}
if $\Lambda/M$ is a free group
(or equivalently if $M=(M\otimes_{\mathbf Z}{\mathbf R})\cap \Lambda$).

We leave the proof of the following well-known result to the reader.

\begin{prop} \label{proporthogreseaux} 
(cf. Chapter I, Proposition 9.8 in \cite{M})
Let $\Lambda$ be an integral
lattice of determinant $1$. 
Let $M,N\subset \Lambda $ be two sublattices of $\Lambda$ such that
$$M=\Lambda\cap N^\perp\hbox{ and }N=\Lambda\cap M^\perp$$
where $X^\perp\subset \Lambda\otimes_{\mathbf Z} {\mathbf R}$ denotes the 
subspace of all vectors orthogonal to $X$ in the Euclidean
vector-space $\Lambda\otimes_{\mathbf Z} {\mathbf R}$ (i.e. $M$ and $N$ are saturated sublattices,
orthogonal to each other and $M\oplus N$ is of finite index in $\Lambda$).

Then the two determinant groups
$$M^\sharp/M\hbox{ and }N^\sharp/N$$
are isomorphic. In particular, the determinants of the
lattices $M$ and $N$ are equal.
\end{prop}

Given a ${\mathbf Z}-$basis $b_1,\dots,b_n\in \Lambda$ of an $n-$dimensional
lattice $\Lambda$, the symmetric positive definite matrix $G$ with coefficients
$$G_{i,j}=\langle b_i,b_j\rangle$$
is called a {\it Gram matrix} of $\Lambda$. Its determinant $\det(G)$
is independent of the choice of the basis and equals the determinant 
of $\Lambda$.

Two lattices $\Lambda$ and $M$ are {\it similar}, 
if there exists a bijection 
$\Lambda\longrightarrow M$ which extends to an Euclidean similarity between
$\Lambda\otimes_{\mathbf Z}{\mathbf R}$ and $M\otimes_{\mathbf Z}{\mathbf R}$. 
The set of 
similarity classes of lattices is endowed with a natural topology:
a neighbourhood of a lattice $\Lambda$ is given by all
lattices having a Gram matrix in ${\mathbf R}_{>0}\ V(G)$
where $V(G)$ is a neighbourhood of a fixed Gram matrix $G$ of $\Lambda$.

Similar lattices have identical densities and the density function
$\Lambda\longmapsto \Delta(\Lambda)$ is continuous with respect
to the natural topology on similarity classes.

Consider the set ${\cal L}_n$ of all $n-$dimensional lattices of the form
$$\Lambda=\{z\in{\mathbf Z}^{n+1}\ \vert\ \langle z,s\rangle=0\}$$
for $s\in{\mathbf N}^{n+1}\setminus\{0\}$.

\begin{prop} \label{Lestdense} The set ${\cal L}_n$ is dense in the set of
similarity classes of $n-$dimensional Euclidean lattices.
\end{prop}

The upper bound for densities of lattices in ${\cal L}_n$ is thus equal to
the maximum for densities of all $n-$dimensional lattices.

{\bf Proof of Proposition \ref{Lestdense}} 
Given a Gram matrix $G=\langle b_i,b_j\rangle$ with respect to a 
${\mathbf Z}-$basis $b_1,\dots,b_n\in \Lambda$ of an $n-$dimensional lattice 
$\Lambda$, Gram-Schmidt orthogonalization shows that
$$G=L\ L^t$$
where $L=(l_{i,j})_{1\leq i,j\leq n}$ is lower triangular (and invertible).

Choose $\kappa>0$ large and consider the integral lower triangular matrix
$\tilde L(\kappa)$ whose coefficients $\tilde l_{i,j}\in{\mathbf Z}$ 
satisfy
$$\vert \tilde l_{i,j}-\kappa l_{i,j}\vert\leq 1/2$$
and are obtained by rounding off the coefficients of $\kappa L$ 
to the nearest integers.

Define an integral matrix 
$$B(\kappa)=\left(\begin{array}{ccccccc}
\tilde l_{1,1}&1&0&0\\
\tilde l_{2,1}&\tilde l_{2,2}&1&0&0\\
\vdots&&&\ddots\\
\end{array}\right)$$
of size $n\times (n+1)$ with coefficients
$$b_{i,j}=\left\lbrace\begin{array}{ll}
\tilde l_{i,j}&\hbox{if }j\leq i\\
1&\hbox{if }j=i+1\\
0&\hbox{otherwise .}\end{array}\right.$$ 

It is easy to see that the rows of $B(\kappa)$ span a saturated integral
sublattice $\tilde \Lambda(\kappa)$ of dimension $n$ in ${\mathbf Z}^{n+1}$.
The special form of $B(\kappa)$ shows that there exists an integral 
row-vector
$$v(\kappa)=\left(\begin{array}{c}
1\\-\tilde l_{1,1}\\
\tilde l_{1,1}\tilde l_{2,2}-\tilde l_{2,1}\\
\vdots\end{array}\right)\in {\mathbf Z}^{n+1}$$
such that $B(\kappa)v(\kappa)=0$.   
We have thus 
$$\tilde \Lambda(\kappa)=
v(\kappa)^\perp\cap{\mathbf Z}^{n+1}\subset {\mathbf E}^{n+1}\ .$$

Since $\hbox{lim}_{\kappa\rightarrow\infty}\frac{1}{\kappa^2}\ B(\kappa)$
is given by the matrix $L$ with an extra row of zeros appended,
we have
$$\hbox{lim}_{\kappa\rightarrow\infty}\frac{1}{\kappa^2}\ B(\kappa)
(B(\kappa))^t=G$$
and the lattice $\frac{1}{\kappa}\tilde \Lambda(\kappa)$
converges thus to the lattice $\Lambda$ for $\kappa\rightarrow \infty$.
Considering the integral vector $s=(s_0,s_1,\dots)$ defined by
$s_i=\vert v(\kappa)_{i+1}\vert$ for $i=0,\dots,n$, we get an integral
lattice 
$$\{z=(z_0,\dots,z_n)\in{\mathbf Z}^{n+1}\ \vert\ \langle 
z,s\rangle=0\}$$
which is isometric to $\tilde \Lambda(\kappa)$.\hfill$\Box$

\section{Proof of Theorem \ref{mainA}}

\begin{lem} \label{lemZn} The standard Euclidean lattice ${\mathbf Z}^n$ contains at
most
$$2\sqrt{\mu+n/4}^n\frac{\pi^{n/2}}{(n/2)!}$$
vectors of (squared Euclidean) norm $\leq \mu$.
\end{lem}

{\bf Proof} We denote by $$B_{\leq \sqrt \rho}(x)=\{z\in{\mathbf E}^n\ 
\vert\ \langle z-x,z-x\rangle\leq \rho\}$$
the closed Euclidean ball with radius $\sqrt{\rho}\geq 0$
and center $x\in{\mathbf E}^n$. Given $0\leq \sqrt{\mu},\ \sqrt{\rho}$
and $x\in B_{\leq \sqrt\mu}(0)$, the closed half-ball
$$\{z\in{\mathbf E}^n\ \vert\ \langle z,x\rangle\leq \langle x,x\rangle\}
\cap B_{\leq \sqrt\rho}(x)$$
(obtained by intersecting the closed halfspace $H_x=
\{z\in{\mathbf E}^n\ \vert\ \langle z,x\rangle\leq \langle x,x\rangle\}$ with the Euclidean ball $B_{\leq \sqrt \rho}(x)$ centered at $x\in
\partial H_x$) is contained in $B_{\leq \sqrt{\mu+\rho}}(0)$.

Since the regular standard cube
$$C=[-\frac{1}{2},\frac{1}{2}]^n\subset {\mathbf E}^n$$
of volume 1 is contained in a ball of radius $\sqrt{n/4}$ 
centered at the origin, 
the intersection
$$(z+C)\cap\{x\in {\mathbf E}^n\ \vert \ \langle x,x\rangle\leq \mu+n/4\}=
(z+C)\cap B_{\leq \sqrt{\mu+n/4}}(0)$$
is of volume at least $1/2$ for any element $z\in {\mathbf E}^n$ of
norm $\langle z,z\rangle \leq \mu$.

Since integral translates of $C$ tile ${\mathbf E}^n$, we have
$$\frac{1}{2}\sharp\{z\in{\mathbf Z}^n\ \vert\ \langle z,z\rangle\leq \mu\}\leq
\hbox{Vol}\ \{x\in{\mathbf E}^n\ \vert\ \langle x,x\rangle\leq \mu+n/4\}
\ .$$
Using the fact that the unit ball in Euclidean 
$n-$space has volume $\pi^{n/2}/(n/2)!$ (cf. 
Chapter 1, formula 17 in \cite{CS}) we get the result. \hfill $\Box$

{\bf Proof of Theorem \ref{mainA}} For $n=0$, the first inequality boils
down to $s_0=1\leq 1+\sqrt{\mu-2}$ and holds for
$\mu\geq 2$. Consider now for $n\geq 1$ a $\mu-$sequence
$(s_0,\dots,s_{n-1})\in {\mathbf N}^n$.

Introduce the set
$${\cal F}=\{(a,k)\in{\mathbf N}^2\ \vert\ \exists\  z=(z_0,\dots,z_{n-1})\in
{\mathbf Z}^n\setminus\{0\}\hbox{ such that}$$
$$ak=\vert\langle z,(s_0,\dots,s_{n-1})\rangle\vert\hbox{ and }
\langle z,z\rangle+k^2<\mu\}\ .$$

Since $\Lambda_{n-1}$ has minimum $\geq \mu$, the equality
$\langle z,(s_0,\dots,s_{n-1})\rangle=0$ implies $\langle z,z\rangle
\geq \mu$ for $z\in {\mathbf Z}^n\setminus\{0\}$.
This shows that we have $a,k>0$ for $(a,k)\in{\mathcal F}$.

Since for a given pair of opposite non-zero 
vectors $\pm z\in{\mathbf Z}^n$ with norm
$0<\langle z,z\rangle<\mu$ there are at most 
$\sqrt{\mu-1-\langle z,z\rangle}\leq \sqrt{\mu-2}$ strictly positive integers
$k$ such that $\langle z,z\rangle+k^2<\mu$, such a pair $\pm z$
of vectors contributes at most $\sqrt{\mu-2}$
distinct elements to $\cal F$. The cardinality $f=
\sharp({\cal F})$
of $\cal F$ is thus bounded by
$$f\leq \sqrt {\mu-2}\frac{\sharp\{z\in{\mathbf Z}^n\ \vert \ 
0<\langle z,z\rangle\leq \mu-1\}}{2}\leq 
\sqrt {\mu-2}\sqrt{\mu-1+n/4}^n\frac{\pi^{n/2}}{(n/2)!}$$
where the last inequality follows from Lemma 3.1. 
There exists thus a strictly positive integer 
$$s_n\leq f+1\leq 1+\sqrt {\mu-2}\sqrt{\mu-1+n/4}^n\frac{\pi^{n/2}}{(n/2)!}$$
such that $(s_n,k)\not\in {\cal F}$ for all $k\in {\mathbf N}$. 
The strictly positive integer $s_n$
satisfies the first inequality of the Theorem and it is 
straightforward to check that the $n-$dimensional lattice
$$\Lambda_n=\{z\in{\mathbf Z}^{n+1}\ \vert\ \sum_{i=0}^n s_iz_i=0\}$$
has minimum $\geq \mu$. This shows the first inequality.

The second inequality
$$1+\sqrt{\mu-2}\sqrt{\mu-1+n/4}^n\frac{\sqrt{\pi}^n}{(n/2)!}
\leq \sqrt{\mu}\sqrt{\mu+n/4}^n\frac{\sqrt{\pi}^n}{(n/2)!}$$
of Theorem \ref{mainA} boils down to 
$$1\leq \sqrt 2\sqrt{2+n/4}\frac{\sqrt{\pi}^n}{(n/2)!}$$
for $\mu=2$. This inequality is clearly true since the $n-$dimensional
Euclideean ball of radius $\sqrt{2+n/4}$ has volume
$\sqrt{2+n/4}\frac{\sqrt{\pi}^n}{(n/2)!}$ and contains the regular
cube $[-\frac{1}{2},\frac{1}{2}]^n$ of volume 1.

For $\mu\geq 3$ we have to establish the inequality
$\Phi(1)-\Phi(0)\geq 1$ where
$$\Phi(t)=\sqrt{\mu-2+2t}\sqrt{\mu-1+t+n/4}^n\frac{\sqrt{\pi}^n}{
(n/2)!}\ .$$
We get thus
$$\begin{array}{lcl}
\Phi(1)-\Phi(0)&\geq&\hbox{inf}_{\xi\in(0,1)}\Phi'(\xi)\\
&\geq&
\frac{1}{\sqrt\mu}\sqrt{\mu-1+n/4}^n\frac{\sqrt{\pi}^n}{(n/2)!}
+\frac{n}{2}\sqrt{\mu-2}\sqrt{\mu-1+n/4}^{n-2}\frac{\sqrt{\pi}^n}{(n/2)!}\ .
\end{array}$$

For $n=1$ and $\mu\geq 2$ we have 
$$\Phi(1)-\Phi(0)\geq\sqrt{1-\frac{1}{\mu}}\frac{\sqrt{\pi}}{\sqrt{\pi}/2}\geq 
\frac{2}{\sqrt 2}>1\ .$$

For $n\geq 2$ and $\mu\geq 3$ we get
$$\Phi(1)-\Phi(0)\geq \sqrt{2+n/4}^{n-2}\frac{\sqrt{\pi}^{n-2}}{((n-2)/2)!}\ 
\pi$$
and the right-hand side equals $\pi>1$ for $n=2$. 
For $n>2$, the right hand side equals $\pi$ times the volume of the
$(n-2)-$dimensional ball of radius $\sqrt{2+n/4}$ containing
the regular cube $[-\frac{1}{2},\frac{1}{2}]^{n-2}$ of volume $1$.
The second inequality follows.
\hfill $\Box$

{\bf Proof of Corollary \ref{corA}} Theorem \ref{mainA} shows the
existence of a $\mu-$sequence $(s_0=1,\dots,s_n)$ satisfying
$$s_0,\dots,s_n\leq \sqrt\mu\sqrt{\mu+n/4}^n\frac{{\sqrt \pi}^n}{(n/2)!}\ .$$
This shows for the lattice 
$\Lambda_n=(s_0,\dots,s_n)^\perp\cap{\mathbf Z}^{n+1}$ the inequality
$$\hbox{det }\Lambda_n=
\sum_{i=0}^n s_i^2\leq (n+1)\mu(\mu+n/4)^n\frac{\pi^n}{((n/2)!)^2}=
(n+1)\mu(\mu+n/4)^n V_n^2$$
and implies
$$\Delta(\Lambda_n)\geq\sqrt{\frac{\mu^n}{4^n(n+1)\mu(\mu+n/4)^n V_n^2}}V_n
\ .$$
This proves the Corollary.\hfill $\Box$

\section{Proof of Theorem \ref{mainB}}

The main idea for proving Theorem \ref{mainB} is to get rid of a factor 
$\sqrt\mu$ when computing an upper bound $f$ for
the size of the finite set $\mathcal F$ considered above during the
proof of Theorem \ref{mainA}. This is possible since the volume 
of the $n-$dimensional unit-ball centered at the origin is concentrated along 
hyperplanes for large $n$. 
For the sake of simplicity, 
we consider sequences in the $\mu\rightarrow\infty$
limit. This allows us to neglect boundary effects when replacing
counting arguments by volume-computations.

In the sequel we write
$$g(x)\sim_{x\rightarrow\alpha}h(x)\ ,\hbox{ respectively }g(x)\leq_{x\rightarrow\alpha}h(x)\ ,$$
for $$\hbox{lim}_{x\rightarrow\alpha}\frac{g(x)}{h(x)}=1\ ,\hbox{ 
respectively }\hbox{limsup}_{x\rightarrow\alpha}\frac{g(x)}{h(x)}\leq 1\ ,$$
where $g(x),\ h(x)>0$.

The following easy Lemma will be useful.

\begin{lem} \label{Vn/Vn-1} We have
$$\sqrt n\ \frac{V_n}{V_{n-1}}=\sqrt{2\pi}(1-\frac{1}{4n}+O(\frac{1}{n^2}))
\ .$$
\end{lem}

{\bf Proof.} Using the definition $V_n=\frac{\pi^{n/2}}{
(n/2)!}$ and Stirlings formula $n!\sim \sqrt{2\pi n}\ n^n\ e^{-n}(1+\frac{1}
{12n}+O(\frac{1}{n^2}))$, we have
$$\begin{array}{lcl}
\sqrt n\ \frac{V_n}{V_{n-1}}&=&
\sqrt{n}\frac{\pi^{n/2}}{\pi^{(n-1)/2}}\ \frac{((n-1)/2)!}{(n/2)!}\\
&=&\sqrt{n\ \pi}\frac{\sqrt{2\pi(n-1)/2}}{\sqrt{2\pi n/2}}\ \frac{
(n-1)^{(n-1)/2}}{2^{(n-1)/2}}\ \frac{2^{n/2}}{n^{n/2}}\ 
\frac{e^{n/2}}{e^{(n-1)/2}}\\
& &\quad \left(\frac{(1+\frac{1}{6(n-1)}+O(\frac{1}{(n-1)^2}))}{
(1+\frac{1}{6n}+O(\frac{1}{n^2}))}\right)\\
&=&\sqrt{2\pi e}(1-\frac{1}{n})^{n/2}(1+O(\frac{1}{n^2}))\\
&=&\sqrt{2\pi e}\ e^{-\frac{n}{2}(\frac{1}{n}+\frac{1}{2n^2}+O(\frac{1}{n^3}))}
(1+O(\frac{1}{n^2}))\\
&=&\sqrt{2\pi}(1-\frac{1}{4n}+O(\frac{1}{n^2}))
\end{array}$$
which ends the proof.\hfill $\Box$

{\bf Proof of Theorem \ref{mainB}} Let $(s_0,\dots,s_{n-1})$
be a finite $\mu-$sequence. 
%Define $\sigma_0,\dots,\sigma_i=\frac{s_i}{\mu^{i/2}V_i},\dots
%\sigma_{n-1}\in{\mathbf R}$ such that
%$$s_i=\sigma_i\mu^{i/2}V_i=\sigma_i\mu^{i/2}\frac{\pi^{i/2}}{(i/2)!}\ .$$
%The reader should think of the real numbers $\sigma_i$ as beeing of order $O(1)$.
For $\epsilon>0$ fixed and suitable $\sigma_n>0$, we show the existence
of a $\mu-$sequence $(s_0,\dots,s_{n-1},s_n)$
with $s_n\in I\cap {\mathbf N}$ where
$I=[\sigma_n\mu^{n/2}V_n,(1+\epsilon)\sigma_n\mu^{n/2}V_n]$.

For $k=1,2,\dots,\in{\mathbf N}$ we define finite subsets
$$I_k=\{s\in I\cap {\mathbf N} \vert\ \sum_{i=0}^{n-1}s_ix_i=ks 
\hbox{ for some }(x_0,\dots,x_{n-1})
\in B_{<\sqrt{\mu-k^2}}\cap {\mathbf Z}^n\}$$
of natural integers in $I\cap {\mathbf N}$ where
$B_{<\sqrt{\mu-k^2}}\cap {\mathbf Z}^n$ denotes the set of 
all integral vectors $(x_0,\dots,x_{n-1})\in {\mathbf Z}^n$ having
(squared Euclidean) norm strictly smaller than $\mu-k^2$.
A sequence $(s_0,\dots,s_{n-1},s_n)$ with $s_n\in I$
is a $\mu-$sequence if 
and only if $s_n\not\in I_k$ for $k=1,2,\dots,\lfloor\sqrt\mu\rfloor$.
Introducing the sets 
$$X_k=\{(x_0,\dots,x_{n-1})\in{\mathbf Z}^n\ \vert\ 
\frac{1}{k}\sum_{i=0}^{n-1}s_ix_i\in I\cap{\mathbf N},
\ \sum_{i=0}^{n-1} x_i^2<\mu-k^2\}$$
we have
$$\sharp(\cup_{k=0}^{\lfloor\sqrt\mu\rfloor} I_k)\leq
\sum_{k=0}^{\lfloor\sqrt\mu\rfloor} \sharp(I_k)\leq
\sum_{k=0}^{\lfloor\sqrt\mu\rfloor} \sharp(X_k)\ .$$
For $\mu$ large enough this ensures the existence of a $\mu-$sequence 
$(s_0,\dots,s_{n-1},s_n)$
with $s_n\leq (1+\epsilon)\sigma_n\mu^{n/2}V_n$ if 
\begin{equation}
(1+\epsilon')\sum_{k=0}^{\lfloor\sqrt\mu\rfloor} \sharp(X_k) 
\leq_{\mu\rightarrow\infty}\epsilon\sigma_n\mu^{n/2}V_n\label{fundineq}
\end{equation}
for some $0<\epsilon'$.

Set 
$$X_k(*)=\{(x_0,\dots,x_{n-1})\in{\mathbf Z}^n\ \vert\ 
\frac{1}{k}\sum_{i=0}^{n-1}s_ix_i\in I,\ \sum_{i=0}^{n-1} x_i^2<\mu-k^2\}$$
and consider the partition
$X_k(*)=X_k(0)\cup X_k(1)\cup\dots\cup X_k(k-1)$ defined by the 
disjoint subsets
$$X_k(a)=\{\sharp\{(x_0,\dots,x_{n-1})\in X_k(*)\ \vert\ 
\sum_{i=0}^{n-1}s_ix_i\equiv
a\pmod k\}\subset X_k(*)\ .$$
Since $s_0=1$ we have (for $\epsilon>0$ fixed) the asymptotic equalities
$$\sharp(X_k(j))\sim_{\mu\rightarrow\infty} \frac{1}{k}\sharp(X_k(*))$$
for $j=0,\dots,k-1$. 

The obvious identity $X_k=X_k(0)$ yields thus 
$$\begin{array}{rcl}
\displaystyle\sharp(X_k)&\displaystyle\sim_{\mu\rightarrow\infty}&\displaystyle\frac{1}{k} 
\sharp\{(x_0,\dots,x_{n-1})\in{\mathbf Z}^n\ \vert\ 
\frac{1}{k}\sum_{i=0}^{n-1}s_ix_i\in I,\ \sum_{i=0}^{n-1} x_i^2<\mu-k^2\}\\
&\displaystyle\sim_{\mu\rightarrow\infty}&\displaystyle \frac{1}{k}
\hbox{Vol}\{(t_0,\dots,t_{n-1})\in {\mathbf E}^n\ \vert \ 
\sum_{i=0}^{n-1} t_i^2\leq \mu-k^2,\ \frac{1}{k}\sum_{i=0}^{n-1}s_it_i\in I
\}\ .\end{array}$$

Setting $\tilde s_{n-1}=\sqrt{\sum_{i=0}^{n-1}s_i^2}$ and $\tilde \sigma_{n-1}=
\frac{\tilde s_{n-1}}{{\sqrt\mu}^{n-1}V_{n-1}}$ we have
$$\sharp(X_k)\sim_{\mu\rightarrow\infty}\frac{1}{k}
\int_\alpha^{(1+\epsilon)\alpha}\sqrt{\mu-t^2}^{n-1}dtV_{n-1}
\leq\frac{\epsilon\alpha}{k}
\sqrt{\mu-\alpha^2}^{n-1}V_{n-1}$$
where $$\alpha=k\frac{\sigma_n\mu^{n/2}V_n}{\tilde s_{n-1}}
= k\sqrt\mu\frac{\sigma_nV_n}{\tilde\sigma_{n-1} V_{n-1}}\ .$$
% {\tilde \sigma_{n-1}\mu^{(n-1)/2}V_{n-1}}

We have thus
$$\sharp(X_k)\leq_{\mu\rightarrow\infty}\epsilon{\sqrt\mu}^n\frac{\sigma_nV_n}{\tilde \sigma_{n-1}}\sqrt{1-k^2\left(\frac{\sigma_nV_n}{\tilde 
\sigma_{n-1} V_{n-1}}\right)^2}^{n-1}$$
implying
$$\sum_{k=1}\sharp(X_k)\leq_{\mu\rightarrow\infty} 
\epsilon{\sqrt\mu}^n\sigma_nV_n\ \frac{1}{\tilde \sigma_{n-1}}
\sum_{k=1}^\infty
\sqrt{1-k^2\left(\frac{\sigma_nV_n}{\tilde 
\sigma_{n-1} V_{n-1}}\right)^2}^{n-1}$$
and showing that the asymptotic inequality (\ref{fundineq}) above 
is satisfied for all $\epsilon>0$ if
\begin{equation}\sum_{k=1}^\infty
\sqrt{1-k^2\left(\frac{\sigma_nV_n}{\tilde 
\sigma_{n-1} V_{n-1}}\right)^2}^{n-1}\leq\frac{\tilde\sigma_{n-1}}{1+\epsilon'}
\ .\label{fineinequation}\end{equation}

Using Lemma \ref{Vn/Vn-1} we get the asymptotics 
$$\sqrt{1-k^2\left(\frac{\sigma_nV_n}{\tilde \sigma_{n-1} V_{n-1}}\right)^2}^{n-1}=
\sqrt{1-k^2\left(\frac{\sigma_n}{\tilde \sigma_{n-1}}\right)^2\frac{2\pi}{n}(1-\frac{1}{2n}+
O(\frac{1}{n^2}))}^{n-1}$$
$$=e^{-k^2\left(\frac{\sigma_n}{\tilde \sigma_{n-1}}\right)^2\pi}\left(
1+k^2\left(\frac{\sigma_n}{\tilde\sigma_{n-1}}\right)^2\pi(\frac{3}{2}-
k^2\left(\frac{\sigma_n}{\tilde \sigma_{n-1}}\right)^2\pi)\frac{1}{n}+O(\frac{1}{n^2})\right)$$
$$<e^{-k^2\left(\frac{\sigma_n}{\tilde \sigma_{n-1}}\right)^2\pi}$$
for $\sigma_n\geq\frac{\tilde\sigma_{n-1}}{\sqrt 2}>\frac{\tilde
\sigma_{n-1}}{k}\sqrt\frac{3}{2\pi}$ and $n$ large enough.

Notice that $\frac{1}{2^{n-1}\tilde \sigma_{n-1}}$ is a lower bound for
the density of the $(n-1)-$dimensional integral lattice
$$\{(x_0,\dots,x_{n-1})\in{\mathbf Z}^n\ \vert\ \sum_{i=0}^{n-1}s_ix_i=0\}
=(s_0,\dots,s_{n-1})^\perp\cap{\mathbf Z}^n$$
associated to the $\mu-$sequence $(s_0,\dots,s_{n-1})$.

For $\sigma_n$ of order $O(1)$ we have
\begin{equation}s_n\sim_{\mu\rightarrow\infty} \tilde s_n=\sqrt{\sum_{i=0}^n
s_i^2}\ .\label{equivalence}\end{equation}
Supposing $\tilde \sigma_{n-1}=
\sum_{k=1}e^{-k^2\pi}$ the choice $\sigma_n=\tilde \sigma_{n-1}
(1+\tilde \epsilon)$ 
implies thus the asymptotic inequality (\ref{fundineq})
for $n$ large enough and all $\tilde \epsilon>0$.
The asymptotic equality (\ref{equivalence})
implies now easily the result and the argument can be iterated.

In the case $\tilde \sigma_{n-1}>\sum_{k=1}e^{-k^2\pi}$
choose $\epsilon'$ small enough such that inequation (\ref{fineinequation})
holds for some $\sigma_n<\tilde\sigma_{n-1}$. This implies
that the asymptotic inequality (\ref{fundineq}) is valid 
and a closer inspection shows that we can iterate this construction
using a decreasing sequence $\tilde \sigma_{n-1}>\sigma_n\geq
\sigma_{n+1}\geq\dots$ with limit $\sum_{k=1}e^{-k^2\pi}$.
This proves the result in this case.

The remaining case $\tilde \sigma_{n-1}<\sum_{k=1}e^{-k^2\pi}$
can for instance be treated by 
replacing the $\mu-$sequence $(s_0,\dots,s_{n-1})$ with a $\mu-$sequence
of smaller density.\hfill$\Box$

\begin{rem} (i) Theorem \ref{mainB} can be slightly sharpened in 
 a standard way which yields the $\zeta(n)$ factor in the 
best known bounds for the density of the densest lattice packing.

\ \ (ii) The main error during the proof of Theorem 1.5 occurs during the 
majoration $$\sharp(I_k)\leq\sharp(X_k)$$
which is very crude. 

\ \ (iii) Instead of working with sublattices of ${\mathbf Z}^{n+1}$
orthogonal to a given vector $(s_0,\dots,s_n)\in {\mathbf Z}^{n+1}$, it
is possible to consider sublattices ${\mathbf Z}^{n+a}$
which are orthogonal to a set of $a\geq 2$ linearly independent
vectors in ${\mathbf Z}^{n+a}$. One might also replace the 
standard lattice ${\mathbf Z}^{n+1}$ by other lattices, e.g.
sublattices of finite index in ${\mathbf Z}^{n+1}$.

\ \ (iv) Let us conclude by mentioning that extending
finite $\mu-$sequences in an optimal way into longer $\mu-$sequences
amounts geometrically to the familiar process  
of lamination for lattices (see for instance \cite{CS} or \cite{M}). 
The existence of an integer $s\in I\setminus I_1$
implies indeed the existence of a 
point $P\in {\mathbf E}^{n-1}$ which is
far away from any lattice point of the affine lattice 
$\{(x_0,\dots,x_{n-1})\ \vert\sum x_is_i=s\}\subset{\mathbf Z}^n$ 
and corresponds thus to a
\lq\lq hole'' of the lattice.
\end{rem}

I thank J. Martinet, P. Sarnak and J-L. Verger-Gaugry
for helpful comments and interest in  this work.

Roland Bacher, INSTITUT FOURIER, Laboratoire de Math\'ematiques, 
UMR 5582 (UJF-CNRS), BP 74, 38402 St MARTIN  D'H\`ERES Cedex (France),
e-mail: Roland.Bacher@ujf-grenoble.fr

\end{document}